%% file: hellas.tex
\documentclass{amsart}
\input jorge

\newcommand{\Geo}{\mathfrak{Geo}}
\newcommand{\Equ}{\mathfrak{Equ}}
\newcommand{\Curl}{\mathfrak{curl} \,} 
\newcommand{\Sign}{\text{sign}} 
\newcommand{\JCC}{\mathcal{J}}

\begin{document}

\title[]{Geometric Knot Spaces and Polygonal Isotopy}
\author[]{Jorge Alberto Calvo}
\dedicatory{Department of Mathematics \\ Williams College \\
Williamstown, MA 01267 \\ email: jcalvo@williams.edu \\ \vspace{20mm}}


\begin{abstract} The space of $ n $-sided polygons embedded in 
three-space consists of a smooth manifold in which points correspond 
to piecewise linear or ``geometric'' knots, while paths correspond to 
isotopies which preserve the geometric structure of these knots.  The 
topology of these spaces for the case $ n = 6 $ and $ n = 7 $ is 
described.  In both of these cases, each knot space consists of five 
components, but contains only three (when $ n = 6 $) or four (when $ n 
= 7 $) topological knot types.  Therefore ``geometric knot 
equivalence'' is strictly stronger than topological equivalence.  This 
point is demonstrated by the hexagonal trefoils and heptagonal 
figure-eight knots, which, unlike their topological counterparts, are 
not reversible.  Extending these results to the cases $ n \ge 8 $ will 
also be discussed.  \\

\end{abstract}
	
\maketitle	
\vspace{-10mm}	
\section{Introduction}

Consider the sorts of configurations that can be constructed out of a 
sequence of line segments, glued end to end to end to form an embedded 
loop in $ \R^{3}.  $ The line segments might represent bonds between 
atoms in a polymer, segments in the base-pair sequence of a circular 
DNA macromolecule, or simply thin wooden sticks attached with flexible 
rubber joints.  Thus, a spatial polygon of this kind serves as a 
mathematical model for some object which is physically knotted yet 
retains some of the rigidity inherited from the materials from which 
it is built.

\begin{figure}[t]
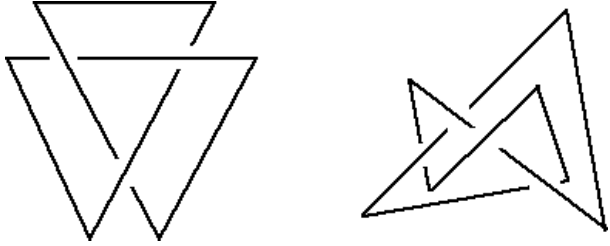

\insertfig{polygons.eps}{1.25in}
\caption{A hexagonal trefoil knot and a heptagonal figure-eight knot.}
\label{fig:polygons}
\end{figure}

 
It is a classical result of three-dimensional topology that knotted 
loops made out of flexible string can always be approximated by 
polygonal loops consisting of many thin, rigid segments.  Furthermore, 
any deformation performed on the string can always be approximated by 
a deformation of the polygon, as long as the number of edges is 
allowed to increase.  However if we insist that the number of edges 
remain constant, then we clearly restrict the types of knots that we 
can construct.  For instance, if we use five or fewer edges, every 
loop we build is topologically unknotted; on the other hand we can 
build a trefoil or a figure-eight knot if we use six or seven edges, 
respectively.  See Figure ~\ref{fig:polygons}.  
What is not clear is whether we can always mimic a topological 
deformation by a deformation of polygons when we place restrictions on 
the number of edges.  For instance, it is unknown whether we can build 
a really complicated polygon which, if it were made out of flexible 
string, could be topologically deformed into a round unknot but, if it 
were built out of rigid sticks with flexible joints, could not be 
flattened out into a planar polygon.  In other words, it is an open 
question whether there exist topological 
unknots which are geometrically knotted.

As it turns out, it is not always possible to find a geometric isotopy 
({\em i.e.} one which keeps the number of edges fixed) between two 
polygonal configurations which are topologically equivalent.  In fact, 
even the case of hexagonal trefoils is nontrivial, as there are 
distinct geometric isotopy types, or {\em isotopes}, of this knot.  As 
a consequence, familiar properties such as reversibility behave 
differently when dealing with geometric knots.


One formulation due to Dick Randell \cite{Randell:conform1, 
Randell:conform2} is obtained by observing the correspondence between 
$ n $-sided polygonal loops in Euclidean three-space and points in $ 
\R^{3n}.  $ Suppose that $ P $ is an $ n $-sided polygon in $ \R^{3}, 
$ together with a choice of a ``first vertex'' $ v_{1} $ and an 
orientation.  By listing the coordinates of each vertex in sequence, 
we obtain a point $ (x_{1}, y_{1}, z_{1}, x_{2}, y_{2}, z_{2}, \ldots, 
x_{n}, y_{n}, z_{n}) \in \R^{3n} $ which we associate with $ P = 
\langle v_{1}, v_{2}, \ldots, v_{n} \rangle.  $ As in the theory of 
Vassiliev invariants, let the {\em discriminant} $ \Sigma^{(n)} $ be 
the set of all points in $ \R^{3n} $ which correspond to polygons with 
self-intersections.  If $ n > 3, $ this discriminant is the union of $ 
\frac{1}{2} n (n - 3) $ pieces, each of which corresponds to the set 
of polygons with an intersecting pair of non-adjacent edges.  For 
instance, the subset in $ \Sigma^{(n)} $ consisting of polygons for 
which the edges $ v_{1}v_{2} $ and $ v_{3}v_{4} $ intersect can be 
described as the collection of polygons for which:
\begin{enumerate} \renewcommand{\labelenumi}{(\roman{enumi})}
\item the vertices $ v_{1}, v_{2}, v_{3}, $ and $ v_{4} $ are coplanar,
\item the line determined by $ v_{1} $ and $ v_{2} $ separates $ v_{3} $ 
from $ v_{4} , $ and
\item the line determined by $ v_{3} $ and $ v_{4} $ separates $ v_{1} $ 
from $ v_{2} . $
\end{enumerate}
Note that this set corresponds to the closure of the locus in $ 
\R^{3n} $ of the system
\begin{gather*}
(v_{2}-v_{1}) \times (v_{3}-v_{1}) \cdot (v_{4}-v_{1}) = 0 , \notag \\
(v_{2}-v_{1}) \times (v_{3}-v_{1}) \cdot (v_{2}-v_{1}) \times 
(v_{4}-v_{1}) < 0 , \\
(v_{4}-v_{3}) \times (v_{1}-v_{3}) \cdot (v_{4}-v_{3}) \times 
(v_{2}-v_{3}) < 0 . \notag
\end{gather*}
Therefore, each of these pieces is the closure of a codimension one 
cubic semi-algebraic variety, {\em i.e.} a hypersurface with boundary.  
We define the space of geometric knots to be the complement of this 
discriminant, $ \Geo^{(n)} = \R^{3n} - \Sigma^{(n)}.  $ Therefore $ 
\Geo^{(n)} $ is a dense open submanifold of $ \R^{3n} .  $ In this 
space, points correspond to embedded polygons or {\em geometric 
knots}, paths correspond to {\em geometric isotopies}, and 
path-components correspond to {\em geometric knot types}.

By a theorem of Whitney \cite{Whitney:varieties}, for any given $ n $ 
there are only finitely many path-components in $ \Geo^{(n)} .$ It is 
also a well-known ``folk theorem,'' due perhaps to Kuiper, that the 
spaces $ \Geo^{(3)}, \Geo^{(4)}, $ and $ \Geo^{(5)} $ are connected.  
In \cite{Calvo:thesis, Calvo:hexagons}, I showed that the spaces $ 
\Geo^{(6)} $ and $ \Geo^{(7)} $ have five components each.  Contrast 
this with the fact that only three topological knot types are 
represented in $ \Geo^{(6)} , $ and that only four topological knot 
types are present in $ \Geo^{(7)} .$ When $ n > 8, $ the exact number 
of path-components remains unknown.  In fact, even the number of 
topological knot types represented in the different components of $ 
\Geo^{(n)} $ is known only when $ n < 9 .$ The following theorem 
summarizes the current status of the classification of geometric knots 
with a small number of edges.

\begin{Thm} \label{thm:classification}
(Calvo \cite{Calvo:thesis, Calvo:hexagons})
\begin{enumerate} \renewcommand{\labelenumi}{(\roman{enumi})}
\item The spaces $ \Geo^{(3)}, \Geo^{(4)}, $ and $ \Geo^{(5)} $ are 
path-connected and consist only of unknots.  
\item The space $ \Geo^{(6)} $ of hexagonal knots contains five 
path-components.  These consist of a single component of unknots, two 
components of right-handed trefoils, and two components of left-handed 
trefoils.
\item The space $ \Geo^{(7)} $ of heptagonal knots contains five 
path-components.  These consist of a single component of unknots and 
of each type of trefoil knot, and two components of figure-eight knots.
\item The space $ \Geo^{(8)} $ of octagonal knots contains at least 
twenty path-components.  However, the only knots represented in this 
space are the unknot, the trefoil knot, the figure-eight knot, every 
five and six crossing prime knot $ (5_1, 5_2, 6_1, 6_2, $ and $ 6_3) , 
$ the square and granny knots $ (3_1 \pm 3_1) , $ the (3, 4)-torus 
knot $ (8_{19}) , $ and the knot $ 8_{20} .  $
\end{enumerate} \end{Thm}

It is important to note that although the deformations obtained as 
paths in $ \Geo^{(n)} $ preserve the polygonal structure of the knot 
in question, in general they will not preserve edge length.  Let $ f: 
\Geo^{(n)} \rightarrow \R^{n} $ be the map taking $ P = \langle v_{1}, 
v_{2}, \ldots, v_{n} \rangle $ to the $ n $-tuple \[ 
(\norm{v_{1}-v_{2}}, \norm{v_{2}-v_{3}}, \ldots, \norm{v_{n-1}-v_{n}}, 
\norm{v_{n}-v_{1}}) .  \] Then points in the preimage $ \Equ^{(n)} = 
f^{-1}(1, 1, \ldots, 1) $ correspond to {\em equilateral knots} with 
unit length edges.  Since the point $ (1, 1, \ldots, 1) $ is a regular 
value for $ f, $ the space $ \Equ^{(n)} $ is a $ 2n $-dimensional 
submanifold (in fact, a codimension $ n $ quadric hypersurface) 
intersecting a number of the components of $ \Geo^{(n)} , $ some 
perhaps more than once.  Paths in this submanifold correspond to 
geometric isotopies which do preserve edge length, so the 
path-components of this space offer yet another notion of knottedness.

In his original papers on molecular conformation spaces 
\cite{Randell:conform1, Randell:conform2}, Randell shows that if $ n 
\le 5 $ then $ \Equ^{(n)} $ is connected.  The case when $ n = 6 $ had 
virtually remained untouched for ten years, except for work by Kenneth 
Millett and Rosa Orellana showing that $ \Equ^{(6)} $ contains a 
single component of topological unknots.~\!\!\footnote{\, Their 
unpublished result is mentioned in Proposition 1.2 of 
\cite{Millett:random.knot}.} By focusing attention to a special case 
of singular ``almost knotted'' hexagons, \cite{Calvo:hexagons} shows 
that two hexagons are equilaterally equivalent exactly when they are 
geometrically equivalent.  Thus $ \Equ^{(6)} $ intersects each 
component of $ \Geo^{(6)} $ exactly once.  Furthermore, 
\cite{Calvo:hexagons} shows that this correspondence of 
path-components is not uninteresting, as the inclusion $ \Equ^{(6)} 
\hookrightarrow \Geo^{(6)} $ has a nontrivial kernel at the level of 
fundamental group.  In fact, if $ \mathcal{T} $ is a component of 
trefoils in $ \Geo^{(6)} ,$ then $ \pi_{1}(\mathcal{T}) = \Z_{2} $ 
while $ \pi_{1}(\mathcal{T} \cap \Equ^{(6)}) $ contains an infinite 
cyclic subgroup.

This paper presents two key ingredients from \cite{Calvo:thesis, 
Calvo:hexagons} used to obtain Theorem~\ref{thm:classification}.  In 
Section~\ref{sec:stratify}, we discuss a method of decomposing $ 
\Geo^{(n)} $ into three-dimensional fibres or ``strata.''  This method 
proved particularly useful in the analysis of $ \Geo^{(6)} $ and $ 
\Geo^{(7)} .$ Then Section ~\ref{sec:project} describes an upper bound 
on the minimal crossing number of the knot realized by an $ n $-sided 
polygon.  This bound, which is obtained by looking at a particular 
projection of polygon into a sphere, improves the one previously known 
by a linear term and provides enough control to classify the 
topological knot types present in $ \Geo^{(8)} .  $

\section{A Stratification of Geometric Knot Spaces} \label{sec:stratify}

Consider the map $ g $ with domain $ \Geo^{(n)} $ which ``forgets'' 
the last vertex of a polygon, mapping 
\[ P = \langle v_{1}, v_{2}, v_{3}, \ldots, v_{n-1}, v_{n} \rangle 
\mapsto g(P) = \langle v_{1}, v_{2}, v_{3}, \ldots, v_{n-1} \rangle .\]
Notice that a generic polygon in $ \Geo^{(n)} $ will map to an 
embedded polygon in $ \Geo^{(n-1)} ; $ the only polygons which do not 
are the ones for which some part of the linkage $ v_{1} v_{2} \ldots 
v_{n-1} $ passes through the line segment between $ v_{1} $ and $ 
v_{n-1}, $ and these polygons form a codimension one subset of $ 
\Geo^{(n)} .$ In particular, since $ \Geo^{(n)} $ is a manifold, every 
$ n $-sided polygon can be perturbed by just a tiny amount so that its 
image under $ g $ lies in $ \Geo^{(n-1)} .  $

Suppose that $ Q $ is an $ (n-1) $-sided polygon in $ \Geo^{(n-1)}.  $ 
Then the preimage $ g^{-1}(Q) $ will be a three-dimensional manifold, 
homeomorphic to the set of valid $ n $th vertices for $ Q. $ This 
divides $ \Geo^{(n)} $ into three-dimensional slices or {\em strata}.  
As $ Q $ varies over $ \Geo^{(n-1)}, $ the corresponding 
three-dimensional stratus $ g^{-1}(Q) $ will vary.  By observing how 
these strata change, we can obtain useful information about $ 
\Geo^{(n)} .$

\begin{figure}[t]
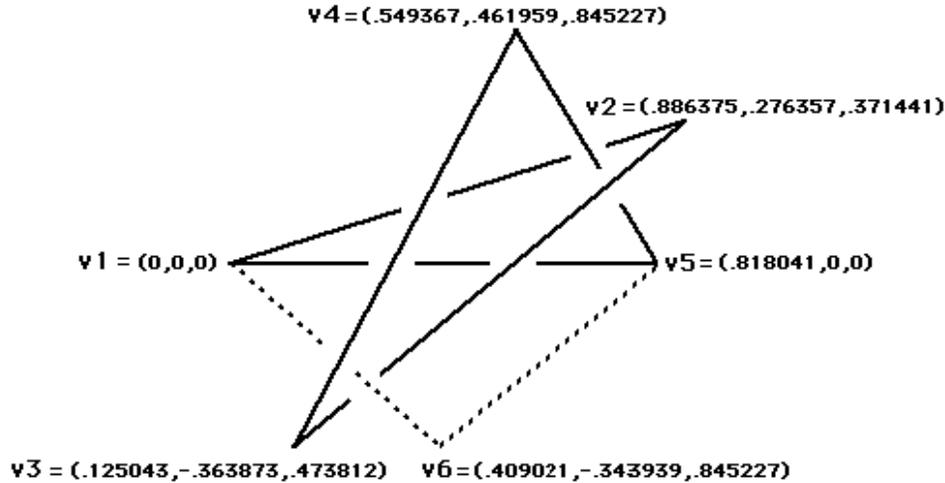
 
\insertfig{frontview.eps}{2.6in} 
\caption{One possible sixth vertex for the pentagon $ Q .$}
\label{fig:frontview} 
\end{figure} 

For example, consider the pentagon $ Q = \langle v_1, v_2, v_3, v_4, 
v_5 \rangle $ with coordinates
\begin{gather*}
\langle (0, 0, 0), (.886375, .276357, .371441), \\
(.125043, -.363873, .473812), \\
(.549367, .461959, .845227), (.818041, 0, 0) \rangle 
\end{gather*}
shown in Figure ~\ref{fig:frontview}.  Suppose that we replace the 
edge between $ v_5 $ and $ v_1 $ with a pair of new edges, from $ v_5 
$ to some new vertex $ v_{6} \in \R^3 $ and from this vertex back to $ 
v_1 .$ This creates a hexagon which, with a bit of care in choosing $ 
v_{6}, $ will also be embedded in $ \R^3 .  $ For instance, if we 
place the new vertex at $ (.4090205, 0, -.912525), $ we obtain an 
unknotted hexagon.  On the other hand, placing $ v_{6} $ at $ 
(.4090205, -.343939, .845227), $ gives a hexagon which is knotted as a 
right-handed trefoil.  See Figure ~\ref{fig:frontview}.  The preimage 
$ g^{-1}(Q) \in \Geo^{(6)} $ is homeomorphic to the dense open subset 
of $ \R^3 $ consisting of ``valid'' sixth vertices for $ Q .$

\begin{figure}[t]
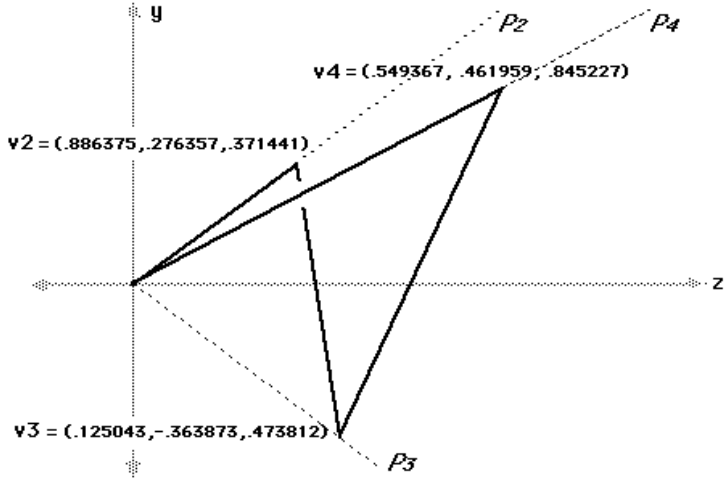
 
\insertfig{sideview.eps}{2.6in}
\caption{Projection  of pentagon $ Q $ into $ yz $-plane.} 
\label{fig:sideview} 
\end{figure}

To examine which points in $ \R^3 $ correspond to embedded hexagons 
obtained from $ Q, $ we will think of the $ x $-axis as a ``central 
axis'' in this space and consider the collection of half-planes 
radiating from this axis.  We refer to these as {\em standard 
half-planes}.  These half-planes appear as rays from the origin in 
Figure ~\ref{fig:sideview}, which shows the projection of $ Q $ into 
the $ yz $-plane.

Let $ \mathcal{P}_{2}, \mathcal{P}_{3}, $ and $ \mathcal{P}_{4} $ be 
the standard half-planes containing $ v_{2}, v_{3}, $ and $ v_{4}, $ 
respectively.  Thus 
\begin{gather*}
\mathcal{P}_{2} = \{ y = \tfrac{276357}{371441} z \approx .744 z, 
\; z > 0 \} , \\
\mathcal{P}_{3} = \{ y = - \tfrac{363873}{473812} z \approx -.768 z, 
\; z > 0 \} , \\
\mathcal{P}_{4} = \{ y = \tfrac{461959}{845227} z \approx .547 z, 
\; z > 0 \} ,
\end{gather*}
as shown in Figure ~\ref{fig:sideview}.

Notice that the interior of any standard half-plane to the left of $ 
\mathcal{P}_{2} $ and $ \mathcal{P}_{3} $ will miss $ Q $ altogether.  
Thus, any point in the interior of any these half-planes may be used 
as a sixth vertex for a hexagon.  Every other standard half-plane, 
however, does intersect $ Q $ at one or more interior points, so these 
half-planes will contain some points which correspond to hexagons with 
self-intersections.

\begin{figure}[t]
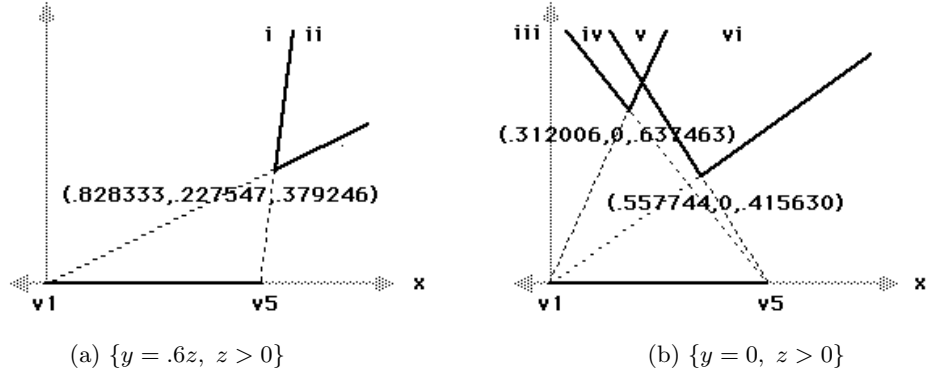
 
\insertfig{vdiscrim.eps}{1.7in}
\centering{{\small (a) $ \{ y = .6z, \; z > 0 \} $ \hspace{46mm} 
(b) $ \{ y = 0, \; z > 0 \} $ \hspace{4mm}}}
\caption{$ Q $ separates each half-plane by ``V''-shaped discriminants.} 
\label{fig:v.discrim} 
\end{figure}

The interior of any standard half-plane between $ \mathcal{P}_{2} $ 
and $ \mathcal{P}_{4} $ will intersect $ Q $ only once, in its second 
edge.  Depending on which point of this plane we choose for the new 
vertex $ v_{6} , $ the two-edge linkage $ v_{5} v_{6} v_{1} $ will 
either dip underneath or jump over this edge.  If $ v_{5} v_{6} v_{1} 
$ goes under $ v_{2} v_{3}, $ then $ v_{6} $ can be dragged back to 
the $ x $-axis, say to the midpoint of edge $ v_{1}v_{5}, $ giving an 
isotopy of the resulting hexagon back to the unknotted loop realized 
by the pentagon $ Q .$ However, if $ v_{5} v_{6} v_{1} $ loops above 
the edge $ v_{2}v_{3} , $ then this edge will obstruct any isotopy of 
the hexagon which attempts to push $ v_{6} $ down towards the $ x 
$-axis in this plane.  For instance, $ Q $ crosses the half-plane $ \{ 
y = .6z, \; z > 0 \} $ at the point $ (.828333, .227547, .379246) .  $ 
Vertices collinear with $ (.828333, .227547, .379246) $ and $ v_{1} $ 
correspond to embedded hexagons only when they lie between these 
points; otherwise the second and sixth edges of the resulting hexagon 
will cross each other.  Similarly, vertices collinear with $ (.828333, 
.227547, .379246) $ and $ v_{5} $ which do not lie between these two 
points correspond to hexagons with intersecting second and fifth 
edges.  Therefore, points in the rays beginning at $ (.828333, 
.227547, $ $ .379246) $ and radiating away from either $ v_{1} $ or $ 
v_{5} $ do not correspond to embedded hexagons, and the half-plane is 
cut into two regions by a ``V''-shaped discriminant.  See Figure 
~\ref{fig:v.discrim}(a).  If $ v_{6} $ is placed in the region of this 
half-plane labelled {\em i}, then the pair of new edges will dip under 
the edge $ v_{2}v_{3} $ and the resulting hexagon will be isotopic to 
$ Q .$ Alternatively, if $ v_{6} $ is placed in the region labelled 
{\em ii}, then $ v_{5} v_{6} v_{1} $ will jump over this edge.

Now, the interior of every standard half-plane between $ 
\mathcal{P}_{4} $ and $ \mathcal{P}_{3} $ intersects $ Q $ in two 
points, in the interior of its second and third edges.  As before, 
these edges will form obstructions to a homotopy moving $ v_{6} $ in 
this plane.  Therefore, for each of the points through which these 
edges cross the half-plane, there will be a ``V''-shaped discriminant 
as above.  For example, in the half-plane $ \{ y = 0, \; z > 0 \} , $ 
which intersects $ Q $ at the points $ (.557744, 0, .415630) $ and $ 
(.312006, 0, .637463) , $ vertices in the four rays beginning at 
either of these points and radiating away from $ v_{1} $ and $ v_{5} $ 
correspond to hexagons with self-intersections.  These two 
``V''-shaped discriminants separate the half-plane into four regions, 
arranged as in Figure ~\ref{fig:v.discrim}(b).  As before, placing the 
new vertex in each of these regions corresponds to looping the new 
edges of the hexagon over either the second ({\em vi}) or third ({\em 
iv}) edge of $ Q , $ or both ({\em v}), or neither ({\em iii}) of 
these.

We can show that the arrangement of the ``V''-shaped discriminants 
remains relatively unchanged for standard half-planes in each of these 
intervals.  In fact, the connected components of the half-planes in 
Figure ~\ref{fig:v.discrim} are only cross-sectional slices of 
``cylindrical sectors'' of $ g^{-1}(Q) $ which wrap around the $ x 
$-axis.  Denote these sectors as {\em i}, {\em ii}, {\em iii}, {\em 
iv}, {\em v}, and {\em vi}, using the notation in Figure 
~\ref{fig:v.discrim}.  Furthermore, let {\em o} denote the sector of $ 
g^{-1}(Q) $ corresponding to vertices in half-planes which do not 
intersect $ Q $ at all.  Then the way in which these sectors are glued 
together depends on the behavior of the discriminants at the three 
``critical level'' standard half-planes $ \mathcal{P}_{2} , 
\mathcal{P}_{3} , $ and $ \mathcal{P}_{4} .$

\begin{figure}[b]
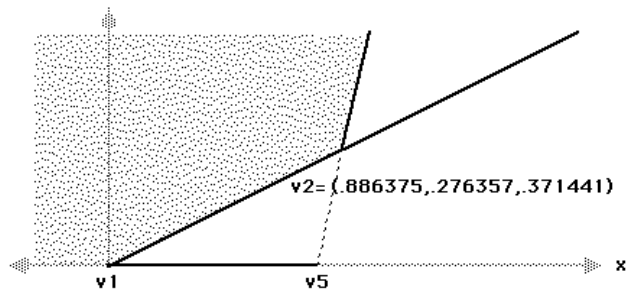
 
\insertfig{p2discrim.eps}{1.54in} 
\caption{Critical level $ \mathcal{P}_{2} = \{ y = 
\tfrac{276357}{371441} z \approx .744 z, \; z > 0 \}.$ }
\label{fig:chp.discrim.a}
\end{figure}

The first of these half-planes, $ \mathcal{P}_{2}, $ contains the 
first edge of $ Q , $ which connects $ v_{1} = (0,0,0) $ to $ v_{2} = 
(.886375, .276357, .371441) .  $ Vertices in rays beginning at any 
point in this edge and radiating away from $ v_{6} $ correspond to 
hexagons with intersecting first and fifth edges.  Hence, for each 
point in this edge there is a ``V''-shaped discriminant.  The union of 
these discriminants forms a two-dimensional discriminant corresponding 
to an obstruction in the space $ g^{-1}(Q) .$ See Figure 
~\ref{fig:chp.discrim.a}.  However, this obstruction only 
partially blocks access to sector {\em i}.  Therefore both {\em i} and 
{\em ii} are glued to {\em o} at this half-plane.

\begin{figure}[h]
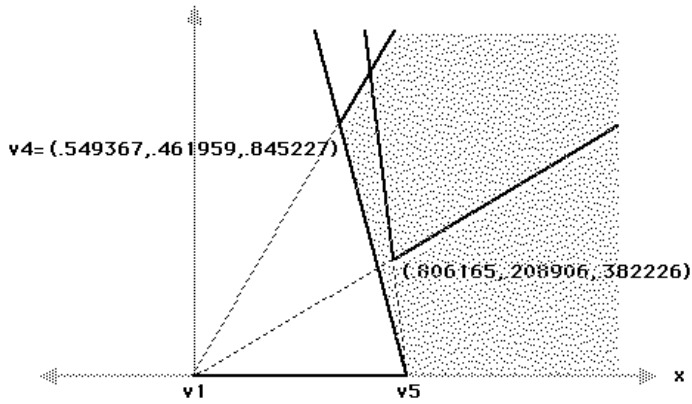
 
\insertfig{p4discrim.eps}{2.13in}
\caption{Critical level $ \mathcal{P}_{4} = \{ y = 
\tfrac{461959}{845227} z \approx .547 z, \; z > 0 \}.$} 
\label{fig:chp.discrim.b} 
\end{figure}

A similar two-dimensional discriminant occurs for $ \mathcal{P}_{4} .$ 
This half-plane contains the fourth edge of $ Q , $ which joins $ 
v_{4} = (.549367, .461959, .845227) $ and $ v_{5} = (.818041, $ $ 0, 
0).$ Vertices collinear with the origin and any point $ p $ on this 
edge correspond to embedded hexagons only if they lie between $ (0, 0, 
0) $ and $ p .$ See Figure ~\ref{fig:chp.discrim.b}.  This 
discriminant completely closes off sector {\em vi}, and obstructs 
parts of sectors {\em i}, {\em ii}, and {\em iii}.  Thus, at this 
level, {\em i} is attached to {\em iii} and {\em iv}, {\em ii} is 
attached to {\em v}, and {\em vi} is abruptly terminated.

The third critical level half-plane, $ \mathcal{P}_{3}, $ presents a 
different situation, as it intersects $ Q $ only at the vertex $ v_{3} 
= (.125043, -.363873, .473812) .$ In this case, the ``V''-shaped 
discriminants corresponding to the second and third edges of $ Q $ 
come together as the two edges become incident at their common vertex.  
As the discriminants merge, sectors {\em iv} and {\em vi} are 
terminated, while both of the sectors {\em iii} and {\em v} merge with 
sector {\em o}.

\begin{figure}[t]
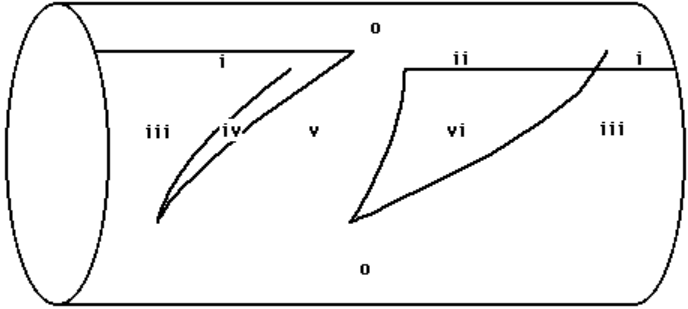
 
\insertfig{cylinder.eps}{2in}
\caption{Cylindrical section of $ \R^{3} $ showing the different 
interconnecting sectors of $ g^{-1}(Q) .$} 
\label{fig:cylinder} 
\end{figure}

Figure ~\ref{fig:cylinder} presents a cylindrical section of $ \R^3 $ 
about the $ x $-axis, showing the sectors of $ g^{-1}(Q) $ and the 
connections between them.  In particular, it shows that $ g^{-1}(Q) $ 
consists of two disjoint path-components, corresponding to the two 
knot types possible for hexagons in the stratus $ g^{-1}(Q) $: the 
unknot and the right-handed trefoil.

The key feature characterizing a pentagon's corresponding stratus in $ 
\Geo^{(6)} $ is the relative position of the second, third, and fourth 
vertices with respect to the axis through the other two vertices.  
Suppose that $ Q = \langle v_1, v_2, v_3, v_4, v_5 \rangle $ is an 
arbitrary pentagon in $ \Geo^{(5)}, $ and that $ \mathcal{L} $ is the 
line determined by $ v_1 $ and $ v_5 .$ Since $ \Geo^{(5)} $ is a 
manifold, we can perturb $ Q $ slightly, if necessary, to ensure that 
$ v_1 v_5 $ is the only edge of $ Q $ which intersects $ \mathcal{L}.$ 
As above, let $ \mathcal{P}_{2}, \mathcal{P}_{3}, $ and $ 
\mathcal{P}_{4} $ be the half-planes with boundary $ \mathcal{L} $ 
which contain $ v_2, v_3, $ and $ v_4 , $ respectively.  Again, a 
slight deformation will make $ Q $ a generic pentagon, guaranteeing 
that the three $ \mathcal{P}_{i} $'s are distinct.

As in the example above, the $ \mathcal{P}_{i} $'s will divide $ 
\R^{3} $ into three open regions, with $ Q $ intersecting two of these 
and completely missing the third.  As we rotate in a right-handed 
fashion about the axis $ \mathcal{L} , $ beginning in the region which 
misses $ Q , $ we will encounter each of the $ \mathcal{P}_{i} $'s in 
one of six orders.  For example, in the pentagon shown in Figures 
~\ref{fig:frontview} and \ref{fig:sideview}, these half-planes appear 
in the order $ \mathcal{P}_2 - \mathcal{P}_4 - \mathcal{P}_3 , $ or 
simply, 2-4-3.

Let $ H \in \Geo^{(6)} $ be a generic hexagon embedded in $ \R^{3} .$ 
By considering the order in which the $ \mathcal{P}_{i} $'s associated 
with $ g(H) $ occur, we divide $ \Geo^{(6)} $ into six open regions 
meeting along codimension one sets where either:
\begin{enumerate} \renewcommand{\labelenumi}{(\roman{enumi})}
\item two of the $ \mathcal{P}_{i} $'s coincide, or
\item an edge of $ H $ crosses $ \mathcal{L} .$
\end{enumerate}
By analyzing the behavior of the strata $ g^{-1}(g(H)) $ as they 
interconnect, we obtain Table ~\ref{tab:hexagons}, which indicates the 
number of path-components in each of the six regions of $ \Geo^{(6)} , 
$ arranged by the topological knot type they represent.

\begin{table}[ht]
\caption{Number of components in each region of $ \Geo^{(6)} .$ }
\label{tab:hexagons}
\begin{center}
\begin{tabular}{|ccccc|}
\hline
$ \mathbf{Region \, of \, \Geo^{(6)}} $ & \; & 
$ \quad \mathbf{0} \quad $ & 
$ \quad \mathbf{3_1} \quad $ & 
$ \; \, \mathbf{- 3_1} \; \, $ \\
\hline 
2-3-4 && 1 & - & - \\
2-4-3 && 1 & 1 & - \\
3-2-4 && 1 & 1 & - \\
3-4-2 && 1 & - & 1 \\
4-2-3 && 1 & - & 1 \\
4-3-2 && 1 & - & - \\
\hline 
\end{tabular}
\end{center}
\end{table}

\begin{figure}[t]
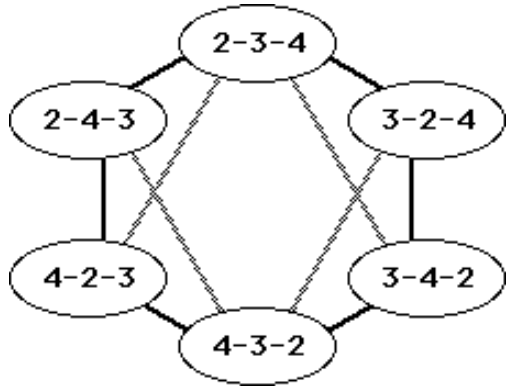

\insertfig{connect.eps}{2in}
\caption{Codimension one connections between regions.}
\label{fig:connections}
\end{figure}

As noted above, the six regions of $ \Geo^{(6)} $ meet along 
codimension one subsets consisting of hexagons for which two of the $ 
\mathcal{P}_i $'s coincide.  For instance, regions 2-4-3 and 4-2-3 
meet along a subset consisting of hexagons with $ \mathcal{P}_2 = 
\mathcal{P}_4 .$ The six regions also meet along codimension one 
subsets consisting of hexagons that intersect line $ \mathcal{L} .  $ 
For example, regions 2-4-3 and 4-3-2 meet along a set of hexagons for 
which edge $ v_2 v_3 $ intersects this line.  These connections are 
shown schematically in Figure~\ref{fig:connections}; here solid lines 
represent hexagons with two coinciding $ \mathcal{P}_i $'s while gray 
lines represent hexagons for which some edge intersects line $ v_1 v_5 
.$

Consider a hexagon $ H $ in the common boundary between two regions of 
$ \Geo^{(6)} .$ Since $ H $ can be perturbed slightly to make generic 
hexagons of either type, $ H $ must be of a topological knot type 
common to both regions.  However, the only knot type common to 
adjacent regions in Figure~\ref{fig:connections} is the unknot.  
Therefore hexagons in these codimension one subsets must be unknotted 
and, in particular, the topological unknots form a single component of 
geometric unknots in $ \Geo^{(6)} .$

On the other hand, suppose that $ h: [0, 1] \rightarrow \Geo^{(6)} $ 
is a path from some trefoil of type 2-4-3 to some trefoil of type 
3-2-4.  Since $ \Geo^{(6)} $ is an open subset of $ \R^{18} , $ there 
is a small open 18-ball contained in $ \Geo^{(6)} $ about each point 
in this path.  Thus we can assume that whenever $ h $ passes through a 
boundary of one of the six regions, it does so through a generic point 
in one of the codimension one subsets above.  But then $ h $ must pass 
through either 2-3-4 or 4-3-2; see Figure~\ref{fig:connections}.  This 
is a contradiction since only unknots live in these regions.  Thus 
there is no path connecting the trefoils of type 2-4-3 and those of 
type 3-2-4.  Similarly, there is no path between the type 4-2-3 and 
3-4-2 trefoils.  This proves that $ \Geo^{(6)} $ consists of five 
path-components: one consisting of unknots, two of right-handed 
trefoils, and two of left-handed trefoils.

The geometric knot types in $ \Geo^{(6)} $ are completely 
characterized by a pair of combinatorial invariants which capture a 
hexagon's topological chirality ({\em i.e.} right- or left-handedness) 
and geometric curl ({\em i.e.} ``upward'' or ``downward'' twisting), 
and are easily computed from the coordinates of a hexagon's vertices.  
To define these invariants, let $ H = \langle v_{1}, v_{2}, v_{3}, 
v_{4}, v_{5}, v_{6} \rangle $ be an embedded hexagon in $ \R^{3} ,$ 
and consider the open triangular disc determined by vertices $ v_{1}, 
v_{2}, $ and $ v_{3}.$ This disc inherits an orientation from $ H $ 
via the ``right hand rule.''  Let $ \Delta_{2} $ be the algebraic 
intersection number of the hexagon and this triangle.  Notice that the 
triangular disc can only be pierced by edges $ v_{4}v_{5} $ and $ 
v_{5}v_{6}.$ Furthermore, if both of these edges intersect the disc, 
they will do so in opposite directions, with their contributions to $ 
\Delta_{2} $ canceling out.  Thus $ \Delta_{2} $ takes on a value of 
0, 1, or $ -1 .$ Similarly, define $ \Delta_{4} $ and $ \Delta_{6} $ 
to be the intersection numbers of H with the triangles $ \triangle 
v_{3}v_{4}v_{5} $ and $ \triangle v_{5}v_{6}v_{1}, $ respectively.  By 
considering the possible values for the $ \Delta_{i} $'s (see Lemma 8 
in \cite{Calvo:hexagons}), we can show that
\begin {enumerate} \renewcommand {\labelenumi} {(\roman{enumi})} 
\item $ H $ is a right-handed trefoil if and only if $ \Delta_2 = 
\Delta_4 = \Delta_6 = 1 , $ 
\item $ H $ is a left-handed trefoil if and only if $ \Delta_2 = 
\Delta_4 = \Delta_6 = -1 , $ and
\item $ H $ is an unknot if and only if $ \Delta_i = 0 $ for some $ i 
\in \{2, 4, 6 \} ,$
\end{enumerate}
implying that the product 
\begin{equation} \label{eq:chirality}
\Delta(H) = \Delta_{2}\Delta_{4}\Delta_{6}, 
\end{equation} 
which we call the {\em chirality} of $ H $, is an invariant under 
geometric deformations.

Next, we define the {\em curl} of $ H $ as
\begin{equation} \label{eq:curl} 
\Curl H = \Sign \bigl( (v_3 - v_1) \times (v_5 - v_1) \cdot (v_2 - 
v_1) \bigr) .
\end{equation} 
This gives the sign of the $ z $-coordinate of $ v_2 $ when we rotate 
$ H $ so that $ v_1, v_3, $ and $ v_5 $ are placed on the $ xy $-plane 
in a counterclockwise fashion, and therefore measures in some sense 
whether a hexagon twists up or down.  Consider a path $ h: [0, 1] 
\rightarrow \Geo^{(6)} $ which changes the curl of a hexagonal trefoil 
from +1 to -1.  Then there must be some point on this path for which 
the vector triple product in (\ref{eq:curl}) is equal to zero.  At 
this point, the vertices $ v_1, v_2, v_3, $ and $ v_5 $ are all 
coplanar.  However, we can show such a hexagon must be unknotted, 
giving us a contradiction.  In particular, the product $ \Delta^{2}(H) 
\, \Curl H $ is also an invariant under geometric deformations.  A 
simple calculation then shows that every trefoil of type 2-4-3 or 
4-2-3 has positive curl, while every one of type 3-2-4 or 3-4-2 has 
negative curl.

\begin{Thm} \label{thm:jcc}
Define the joint chirality-curl of a hexagon $ H $ as the ordered pair 
$ \JCC (H) = (\Delta(H), \Delta^{2}(H) \, \Curl H) .$ Then
\begin{equation} \label{eq:jcc.values}
\JCC (H) = \begin{cases} ( 0, 0 ) & \text{iff} \, H \, \text{is an 
unknot,} \\
( +1, c ) & \text{iff} \, H \, \text{is a right-handed trefoil with} 
\, \Curl H = c, \\
( -1, c ) & \text{iff} \, H \, \text{is a left-handed trefoil with} \, 
\Curl H = c.
\end{cases}
\end{equation}
Therefore the geometric knot type of a hexagon $ H $ is completely 
determined by the value of its chirality and curl.
\end{Thm}

Before leaving the world of hexagons behind, let us make one last 
observation.  Recall that the construction of $ \Geo^{(n)} $ depends 
on a choice of a ``first vertex'' $ v_1 $ and an orientation.  This 
amounts to choosing a sequential labeling $ v_1, v_2, \ldots, v_n $ 
for the vertices of each polygon.  A different choice of labels will 
lead to a different point in $ \Geo^{(n)} $ corresponding to the same 
underlying polygon.  Thus the dihedral group $ \mathbf{D}_n $ of order 
$ 2n $ acts on $ \Geo^{(n)} $ by shifting or reversing the order of 
these labels, and this action preserves topological knot type.  
Observation of the effects on $ \Curl H $ by the group action of $ 
\mathbf{D}_6 $ on $ \Geo^{(6)} $ reveals that same statement does not 
hold true for geometric knot type.  In particular, if the group action 
is defined by the automorphisms
\begin{align*}
r : \langle v_1, v_2, v_3, v_4, v_5, v_6 \rangle & \mapsto \langle v_1, v_6, 
v_5, v_4, v_3, v_2 \rangle , \\
s : \langle v_1, v_2, v_3, v_4, v_5, v_6 \rangle & \mapsto \langle v_2, v_3, 
v_4, v_5, v_6, v_1 \rangle ,
\end{align*}
then 
\[ \Curl r H = \Curl s H = - \Curl H . \]
This shows that the hexagonal trefoil knot is not reversible: In 
contrast with trefoils in the topological setting, reversing the 
orientation on a hexagonal trefoil yields a different geometric knot.  
Furthermore, shifting the labels over by one vertex also changes the 
the knot type of a trefoil, so that taking quotients under this action 
we can see that the spaces $ \Geo^{(6)} / \!  \prec \!  s \!  \succ $ 
of non-based, oriented hexagons, and $ \Geo^{(6)} / \!  \prec \!  r, s 
\!  \succ $ of non-based, non-oriented hexagons consist of only three 
components each.

\begin{figure}[t]
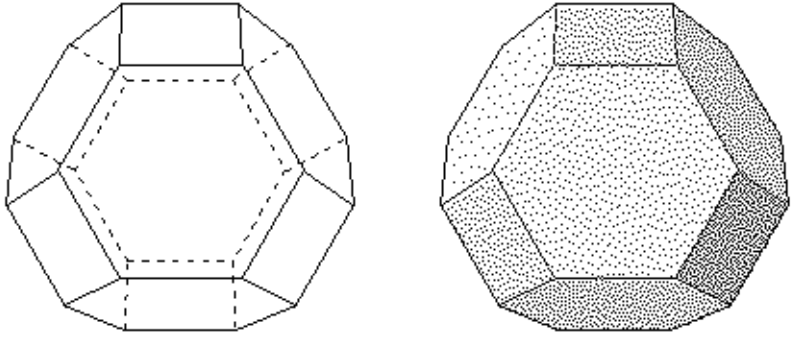

\insertfig{permuta.eps}{2in}
\caption{The valence 3 permutahedron.}
\label{fig:permutahedron}
\end{figure}

A similar decomposition can be made for the space $ \Geo^{(7)} $ of 
heptagons.  In this case we consider the relative ordering of the 
half-planes $ \mathcal{P}_{2}, \mathcal{P}_{3}, \mathcal{P}_{4}, $ and 
$ \mathcal{P}_{5} $ bounded by the line through $ v_{1} $ and $ v_{6}.  
$ This defines 24 open regions which meet along codimension one 
subsets where two of the $ \mathcal{P}_i $'s coincide.  These 
junctions can be schematically described as switches in the indices 
denoting the regions.  For instance, regions 2-4-3-5 and 4-2-3-5 meet 
along a subset consisting of hexagons with $ \mathcal{P}_2 = 
\mathcal{P}_4 .  $ We can build a model for these connections by 
taking a vertex for each of the 24 regions and an edge for each 
codimension-1 subset joining them.  The result is a valence-3 graph 
which forms the 1-dimensional skeleton of a solid zonotope called a 
{\em permutahedron}, shown in Figure ~\ref{fig:permutahedron}.  Each 
vertex of the permutahedron is part of a unique square face 
corresponding to the order-4 sequence of index switches in which the 
first two indices and the last two indices are switched in an 
alternating fashion.  In addition, each vertex is part of two distinct 
hexagonal faces which correspond to the order-6 switch sequences in 
which either the first or last index is fixed while the other three 
indices are permuted through all six possible orderings.  Therefore 
the valence-3 permutahedron has six square faces and eight hexagonal 
faces.  Extending the edges shared by any two hexagonal faces shows 
that this is nothing more than a truncated octahedron, also known in 
crystallography as a Fedorov cubo-octahedron.  ~\!\!\footnote{\, The 
cubo-octahedron is a {\em parallelohedron}, that is, a crystalline 
shape having parallel opposite faces with which three-space can be tiled.  
One should not confuse Fedorov's cubo-octahedron with Kepler's 
cuboctahedron, which is built from an octahedron by truncating at the 
midpoint (rather than at the one- and two-third points) of each edge 
and thus consists of six squares and eight triangular faces.  See 
pp.17 -- 18 in \cite{Ziegler:polytopes} and pp.  722 -- 723 in 
\cite{Tutton:crystals}.}

With a few additional considerations, ~\!\!\footnote{\, The interested 
reader is referred to pp.  53 -- 56 in \cite{Calvo:thesis}.} the 
analysis of the strata over each of these regions shows that $ 
\Geo^{(7)} $ has a single path-component of unknots and of each 
topological type of trefoil, and two containing figure-eight knots.  
Again, these figure-eight knots are new examples of distinct geometric 
isotopes of the same topological knot, which can be distinguished by a 
geometric invariant $ \Xi, $ defined as follows.

Suppose that $ H $ is the heptagon $ \langle v_1, v_2, v_3, v_4, v_5, 
v_6, v_7 \rangle.  $ Define the functions $ \Theta_3 (H) $ and $ 
\Theta_6 (H) $ as
\begin{equation}
\begin{aligned}
\Theta_3 &= \Sign \bigl( (v_7 - v_1) \times (v_2 - v_1) \cdot (v_3 - v_1) \bigr), \\
\Theta_6 &= \Sign \bigl( (v_6 - v_1) \times (v_7 - v_1) \cdot (v_2 - v_1) \bigr).
\end{aligned} 
\end{equation}
Then $ \Theta_3 = \Theta_6 $ if the vertices $ v_3 $ and $ v_6 $ lie 
on the same side of the plane $ \mathcal{P} $ determined by $ v_7, 
v_1, $ and $ v_2 , $ and $ \Theta_3 = - \Theta_6 $ if $ v_3 $ and $ 
v_6 $ lie on different sides of $ \mathcal{P} .  $ Notice that 
for a generic heptagon, exactly one of the functions $ \frac{1}{2} 
(\Theta_3 + \Theta_6) $ and $ \frac{1}{2} (\Theta_3 - \Theta_6) $ is 
zero, while the other is $ \pm 1 .  $

Let $ I_{34} $ denote the algebraic intersection number of edge $ v_3 
v_4 $ with the triangular disc $ \triangle v_{7}v_{1}v_{2}, $ using 
the usual orientations induced by $ H .  $ Similarly define $ I_{45} $ 
and $ I_{56} $ as the intersection numbers of the triangle $ \triangle 
v_{7}v_{1}v_{2} $ with the edges $ v_4 v_5 $ and $ v_5 v_6 , $ 
respectively.

If $ H $ has $ \Theta_3 = \Theta_6 , $ then $ v_3 $ and $ v_6 $ lie on 
the same side of the plane $ \mathcal{P} $ so that the three-edge 
linkage $ v_3 v_4 v_5 v_6 $ will intersect $ \mathcal{P} $ at most 
twice.  Furthermore, if both of these intersections happen in the 
interior of $ \triangle v_{7}v_{1}v_{2}, $ they occur with opposite 
orientations.  Thus the sum $ I_{34} + I_{45} + I_{56} $ only takes on 
values -1, 0, or 1.

\begin{figure}[t]
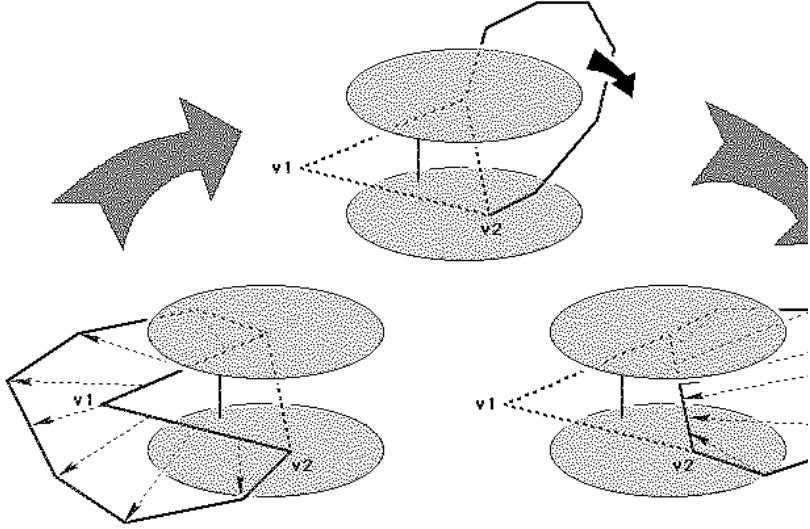

\insertfig{jumprope.eps}{2.75in}
\caption{A piecewise linear isotopy of the linkage $ v_7 v_1 v_2 .$}
\label{fig:jumprope}
\end{figure}

On the other hand, suppose that $ H $ is a figure-eight knot with $ 
\Theta_3 = - \Theta_6 .  $ Then $ v_3 $ and $ v_6 $ lie on the 
opposite sides of the plane $ \mathcal{P} $ so that the three-edge 
linkage $ v_3 v_4 v_5 v_6 $ intersects $ \mathcal{P} $ an odd number 
of times.  First, suppose that there is only one intersection; then 
the linkage $ v_7 v_1 v_2 $ can be piecewise linearly isotoped into a 
straight line segment.  We can think of this isotopy as either pushing 
$ v_{1} $ in a straight line path towards the midpoint of the line 
segment $ v_{2} v_{7}, $ or (in the case that the intersection occurs 
inside $ \triangle v_{7}v_{1}v_{2} $) as stretching $ v_7 v_1 v_2 $ 
into a large loop, swinging it like a ``jump rope'' around and to the 
other side of the heptagon, and then pushing it in until it coincides 
with the line segment $ v_{2} v_{7} .  $ See Figure 
~\ref{fig:jumprope}.  In either case we get a hexagonal realization of 
a figure-eight knot.  Since this is impossible, the linkage $ v_3 v_4 
v_5 v_6 $ has to cross the plane $ \mathcal{P} $ three times, and in 
particular, $ v_3 v_4 $ and $ v_5 v_6 $ must do so with the same 
orientation.  Therefore the quantity $ I_{34} - I_{56} $ will either 
be zero (when both edges intersect $ \triangle v_{7}v_{1}v_{2}, $ or 
when neither of the two do) or $ \pm 1 $ (when only one of these 
intersections occurs inside the triangle).

A quick look at the possible configurations shows that: 
\begin{enumerate} \renewcommand {\labelenumi} {(\roman{enumi})} 
\item if $ H $ is a heptagonal figure-eight knot with $ \Theta_3 = 
\Theta_6 , $ then exactly one of the intersection numbers $ I_{34}, 
I_{45}, $ or $ I_{56} $ is non-zero; thus $ I_{34} + I_{45} + I_{56} = 
\pm 1 $ (Lemma 4.2 in \cite{Calvo:thesis}), and 
\item if $ H $ is a heptagonal figure-eight knot with $ \Theta_3 = - 
\Theta_6 , $ then exactly one of the intersection numbers $ I_{34} $ 
or $ I_{56} $ is non-zero; in particular $ I_{34} - I_{56} = \pm 1 $ 
(Lemma 4.3 in \cite{Calvo:thesis}).
\end{enumerate}
Now consider the function
\begin{equation} \label{eq:xi}
\Xi (H) = \frac{1}{2} \biggl( \Theta_3 + \Theta_6 \biggr) 
\biggr( I_{34} + I_{45} + I_{56} \biggr) +
\frac{1}{2} \biggl( \Theta_3 - \Theta_6 \biggr) 
\biggr( I_{34} - I_{56} \biggr) . 
\end{equation}
By (i) and (ii) above, $ \Xi $ can only take values of 1 or -1.  
Suppose that the value of $ \Xi $ changes along some path $ h: [0,1] 
\rightarrow \Geo^{(7)} .$ Since $ \Geo^{(7)} $ is a manifold, we can 
assume that, in this path, only one vertex passes through the interior 
of $ \triangle v_{7}v_{1}v_{2} $ at any one time, and that only one 
edge intersects the line segment $ v_7 v_2 $ at a time, and that these 
two things happen at different times.  Note that each of these events 
will change the values of $ I_{34} + I_{45} + I_{56} $ and $ I_{34} - 
I_{56} $ by at most one.  However, $ \Xi $ can only change in 
increments of two, so if the values of $ \Theta_3 $ and $ \Theta_6 $ 
remain constant through out $ h, $ $ \Xi $ must also remain unchanged.
 
By reversing orientations if necessary, we can assume then that the 
deformation changes the sign of $ \Theta_3 .  $ In particular, let $ 
H_0 $ be a heptagonal figure-eight knot with $ \Theta_3 = 0 .  $ By 
pushing $ v_3 $ slightly towards $ v_6 , $ we get a heptagon $ H_0^+ $ 
with $ \Theta_3 = \Theta_6 ; $ let $ I_{34}^+ $ be the appropriate 
intersection number for this heptagon.  On the other hand, we obtain a 
heptagon $ H_0^- $ with $ \Theta_3 = - \Theta_6 $ by pushing $ v_3 $ 
away from $ v_6 ; $ let $ I_{34}^- $ be the corresponding intersection 
number for this heptagon.  By picking $ H_{0}^{+} $ and $ H_{0}^{-} $ 
close enough to $ H_{0}, $ we can assume that the values of the 
intersection numbers $ I_{45} $ and $ I_{56} $ coincide for all three 
knots.  This leaves two cases to consider.

First, suppose that $ I_{34}^- = 0 .  $ Then $ I_{56} = \pm 1 $ by 
(ii), and hence $ I_{34}^+ = I_{45} = 0 $ by (i).  Therefore
\[ \biggl( I_{34}^+ + I_{45} + I_{56} \biggr) = I_{56} = 
- \biggl( I_{34}^- - I_{56} \biggr) . \] 
The extra negative sign in the right hand term of this equation 
neutralizes the change of sign in $ \Theta_{3}, $ so that $ \Xi $ 
remains unchanged.

Next, suppose that $ I_{34}^- = \pm 1, $ in which case $ I_{34}^+ = 0 
.$ Then $ I_{56} = 0 $ by (ii) and $ I_{45} = \pm 1 $ by (i).  
Furthermore, the edges $ v_3 v_4 $ and $ v_4 v_5 $ must intersect the 
interior of $ \triangle v_{7}v_{1}v_{2} $ from opposite directions, so 
$ I_{34}^- = - I_{45} .$ Therefore
\[ \biggl( I_{34}^+ + I_{45} + I_{56} \biggr) = I_{45} = - I_{34}^- =
- \biggl( I_{34}^- - I_{56} \biggr) , \]
so that, as before, $ \Xi $ does not change.  This proves the 
following result.

\begin{Thm}
$ \Xi $ is an invariant of heptagonal figure-eight knots under geometric 
deformations.
\end{Thm}

\begin{figure}[p]
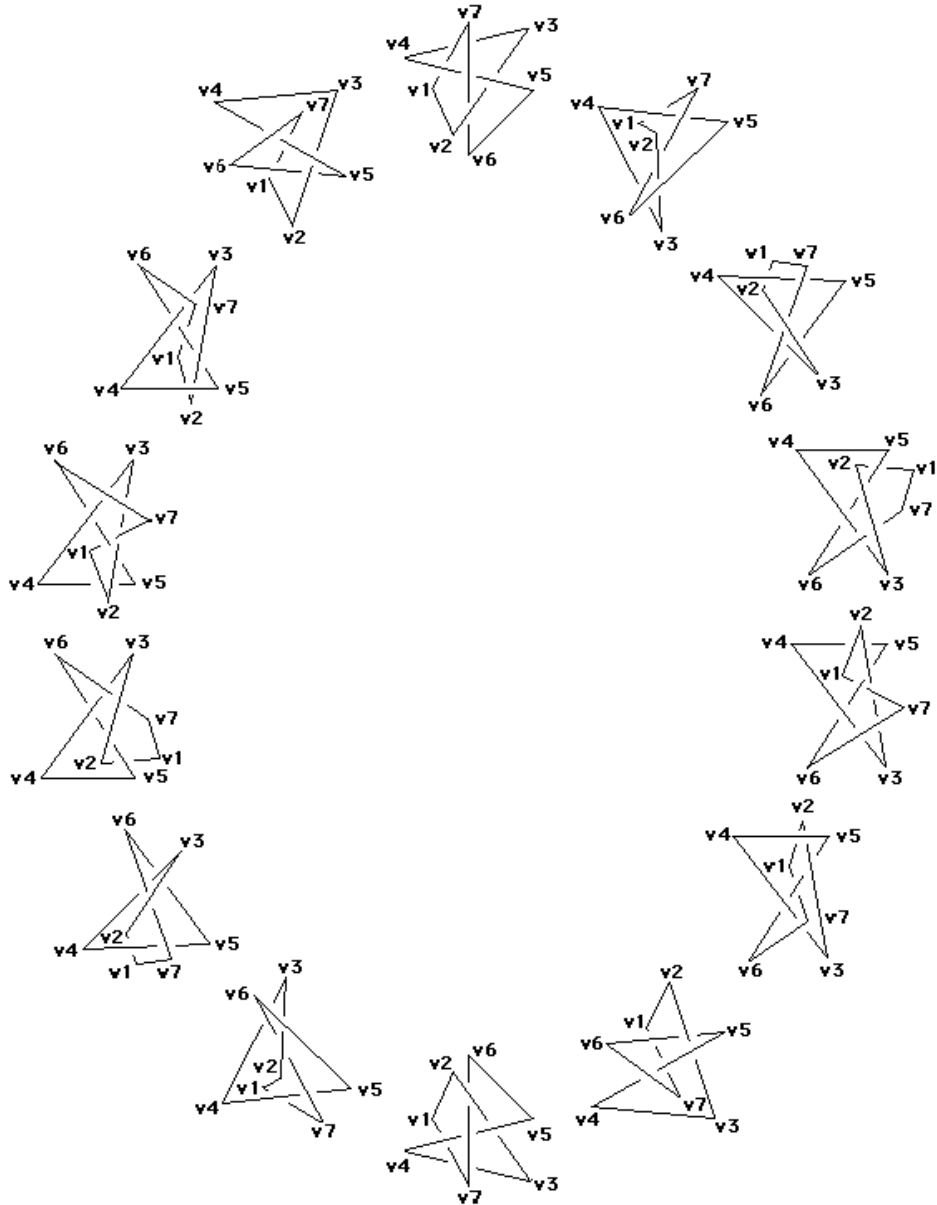
 
\vspace{.5in}
\insertfig{achiral.eps}{6.4in}
\caption{Heptagonal figure-eight knots are achiral.}
\label{fig:achiral}
\end{figure}

It is interesting to note that $ \Xi $ is also invariant under mirror 
reflections, since the resulting sign changes in the functions $ 
\Theta_{3}, \Theta_{6}, I_{34}, I_{45}, $ and $ I_{56} $ cancel out in 
(\ref{eq:xi}).  This reflects the fact that heptagonal figure-eight 
knots are achiral, {\em i.e.} equivalent to their mirror images.  
Figure ~\ref{fig:achiral} shows one such isotopy.  Starting with the 
diagram at the top of Figure ~\ref{fig:achiral} and proceeding in a 
clockwise fashion, we first push $ v_1 $ through the interior of the 
triangular disc $ \triangle v_2 v_3 v_4 .$ Note that in doing so, we 
may need to change the lengths of one or more of the edges.  Although 
it is difficult to see from the perspective of Figure 
~\ref{fig:achiral}, this motion actually defines an isotopy from the 
heptagon $ \langle v_1, v_2, v_3, v_4, v_5, v_6, v_7 \rangle $ to the 
heptagon $ \langle -v_6, -v_7, -v_1, -v_2, -v_3, -v_4, -v_5 \rangle .$ 
We continue by repeating similar moves, passing $ v_3 $ through $ 
\triangle v_4 v_5 v_6 , $ then $ v_5 $ through $ \triangle v_6 v_7 v_1 
, $ and so on.  After seven steps, when we move $ v_6 $ past $ 
\triangle v_7 v_1 v_2 , $ we arrive at the diagram at the bottom of 
Figure ~\ref{fig:achiral}.  At this point, the figure-eight knot is 
the mirror image of the starting position.

Finally, consider the $ \mathbf{D}_{7} $ action on $ \Geo^{(7)} $ defined by 
the automorphisms
\begin{align*}
r \langle v_1, v_2, v_3, v_4, v_5, v_6, v_7 \rangle &= 
\langle v_1, v_7, v_6, v_5, v_4, v_3, v_2 \rangle \\
s \langle v_1, v_2, v_3, v_4, v_5, v_6, v_7 \rangle &= 
\langle v_2, v_3, v_4, v_5, v_6, v_7, v_1 \rangle .
\end{align*}
Reversing the orientation on $ H $ via the map $ r $ will reverse the 
orientations on both the edges of $ H $ and the triangular discs that 
they define.  In particular,
\[ I_{34}(rH) = I_{56}(H) \qquad \qquad I_{45}(rH) = I_{45}(H) 
\qquad \qquad I_{56}(rH) = I_{34}(H) .
\]
On the other hand, $ r $ not only switches the roles of $ \Theta_3 $ and $ \Theta_6, $
but also changes their signs:
\begin{align*}
\Theta_3 (rH) &= \Sign \bigl( (v_2 - v_1) \times (v_7 - v_1) \cdot (v_6 - v_1) \bigr), \\
&= - \Sign \bigl( (v_6 - v_1) \times (v_7 - v_1) \cdot (v_2 - v_1) \bigr), \\
&= - \Theta_6 (H), \\
\Theta_6 (rH) &= \Sign \bigl( (v_3 - v_1) \times (v_2 - v_1) \cdot (v_7 - v_1) \bigr), \\
&= - \Sign \bigl( (v_7 - v_1) \times (v_2 - v_1) \cdot (v_3 - v_1) \bigr), \\
&= - \Theta_3 (H). \\
\end{align*}
Therefore
\begin{align*} 
\Xi (rH) &= \frac{1}{2} \biggl( \Theta_3 (rH) + \Theta_6 (rH) \biggr) 
\biggr( I_{34} (rH) + I_{45} (rH) + I_{56} (rH) \biggr) \\ & \hspace{25mm} +
\frac{1}{2} \biggl( \Theta_3 (rH) - \Theta_6 (rH) \biggr) 
\biggr( I_{34} (rH) - I_{56} (rH) \biggr) \\
&= \frac{1}{2} \biggl( - \Theta_6 (H) - \Theta_3 (H) \biggr) 
\biggr( I_{56} (H) + I_{45} (H) + I_{34} (H) \biggr) \\ & \hspace{30mm} +
\frac{1}{2} \biggl( - \Theta_6 (H) + \Theta_3 (H) \biggr) 
\biggr( I_{56} (H) - I_{34} (H) \biggr) \\
&= - \Xi (H) . 
\end{align*}
This shows that, like hexagonal trefoil knots, figure-eight knots in $ 
\Geo^{(7)} $ are irreversible, in contrast with their topological 
counterparts.  However, recall that the irreversibility of trefoils in 
$ \Geo^{(6)} $ depended strongly on our choice of a ``first'' vertex $ 
v_1 .$ In that case, a cyclic permutation of its six vertices would 
change the trefoil's geometric knot type.  This is not the case for 
the figure-eight knots in $ \Geo^{(7)}, $ for consider the group 
action induced by the automorphism $ s $ on the set of geometric 
isotopes of the figure-eight knot.  This is an order 7 action on a two 
element set, and must therefore be trivial.  In other words, we must 
have 
\[ \Xi (sH) = \Xi (H).  \] 
Hence the distinction in the two figure-eight knot types is an effect 
of ``true'' geometric knotting, which goes beyond a simple relabeling 
of the vertices or our arbitrary choice of first vertex.

\section{Knot Projections and Minimal Polygon index} \label{sec:project}
 
In Section ~\ref{sec:stratify}, we were concerned with the question of 
determining, for a given integer $ n, $ the number of path-components 
present in the space $ \Geo^{(n)} $ of $ n $-sided polygons.  In other 
words, ``how many geometric knot types are there for a particular 
value of $ n $?''  However, as $ n $ increases, the space $ \Geo^{(n)} 
$ becomes more and more combinatorially intricate.  As this happens, 
we turn to the question of understanding the number of represented 
topological (rather than geometric) knot types, and in particular, of 
how complicated a knot can be realized by an $ n $-sided polygon.  The 
answer to this question is only known when $ n \le 8 .$ For example, 
we know there are 9-sided polygonal embeddings of every seven crossing 
prime knot $ ( 7_{1}, \ldots , 7_{7}) $ as well as the knots $ 8_{16}, 
8_{17}, 8_{18}, 8_{21}, 9_{40}, 9_{41}, 9_{42}, $ and $ 9_{46}, $ but 
presumably this list could be much bigger, and include some of the 
knots for which we have so far only found 10- or 11-sided 
realizations.~\!\!\footnote{\, See Table 1 in 
\cite{Calvo:montecarlo}.} In this section, we give one of several 
known bounds on the complexity of an $ n $-sided polygon.

Recall that the {\em minimal crossing number} of a knot is the 
smallest number of crossings present in any general position 
projection of the knot into a plane or sphere.  This is the 
conventional measure of a knot's complexity, used in the standard 
notation for knots and links as well as in the knot tables in the 
appendices of \cite{Adams:book}, \cite{Kauffman:book}, 
\cite{Livingston:book}, and \cite{Rolfsen:book}.  We similarly define 
the {\em minimal polygon index} of a knot as the smallest number of 
edges present in any polygonal embedding of the knot.  This invariant, 
which is elsewhere known as the {\em stick number} \cite{Adams:book, 
Adams:sticks, Furstenberg:sticks}, the {\em broken line number} 
\cite{Negami:s.bounds}, or simply the {\em edge number} 
\cite{Meissen:sticks, Randell:sticks}, serves as the corresponding 
measure of complexity for polygonal knots.  These two invariants are 
traditionally related by the following construction.  
~\!\!\footnote{\, This construction appears in Theorem 7 in 
\cite{Negami:s.bounds}, and as Exercise 1.38 in \cite{Adams:book}.}

Let $ P = \langle v_1, v_2, \ldots, v_{n-1}, v_n \rangle \in 
\Geo^{(n)} $ be an $ n $-sided polygon embedded in $ \R^{3} .$ We 
project the points in $ P $ orthogonally onto a plane perpendicular to 
one of its edges, say $ v_{1} v_{2} .$ This amounts to looking at the 
polygon from a viewpoint in which we see the edge $ v_{1}v_{2} $ 
``head on,'' so that the image of $ P $ on our retina ({\em i.e.} the 
plane) is an $ (n-1) $-sided polygon.  An edge in this image cannot 
cross either of its two neighbors, or itself, so each edge will 
intersect at most $ n-4 $ other edges.  Thus for generic polygons in $ 
\Geo^{(n)}, $ this method gives a knot projection with no more than $ 
\tfrac{1}{2} (n-1)(n-4) $ crossings.  This leads to the conclusion 
that if a knot $ K $ has minimal crossing number $ c(K) $ and minimal 
polygon index $ s(K) ,$ then \[ c(K) \le \frac{\bigl( s(K) - 1 \bigr) 
\bigl( s(K) - 4 \bigr)}{2} , \] or equivalently (by completing squares 
and solving for $ s $) \[ s(K) \ge \frac{5 \, + \, \sqrt{9 \, + \, 8 
c(K)}}{2} .\]

For hexagons and heptagons, the bound on crossing number becomes 5 and 
9, respectively, well over the actual values of 3 and 4 obtained in 
Section ~\ref{sec:stratify}.  In fact, the estimated $ \tfrac{1}{2} 
(n-1)(n-4) $ crossings in the image of an $ n $-sided polygon can 
never be achieved when $ n $ is odd.  Here we present an improvement 
on the bounds above.

First suppose that we relabel the vertices of $ P $ in sequence so 
that $ v_{1} $ is a point on the boundary of the convex hull spanned 
by the vertices of $ P .$ Therefore, we can find a plane $ 
\mathcal{P}_{1} $ which intersects $ P $ only at the vertex $ v_{1}, $ 
with $ P $ lying entirely on one side of $ \mathcal{P}_{1} .$

Let $ \mathcal{S} $ be a large sphere centered at $ v_1 $ and 
enclosing all of $ P , $ and consider the image of the radial 
projection $ p: P - \{v_1\} \rightarrow \mathcal{S} .$ By our choice 
of $ v_{1} ,$ this image lies entirely in a hemisphere of $ 
\mathcal{S} $ cut by the equator $ \mathcal{S} \cap \mathcal{P}_{1} .$ 
Furthermore, note that the interiors of edges $ v_1 v_2 $ and $ v_1 
v_n $ are respectively mapped to the single points $ p(v_2) $ and $ 
p(v_n) .$ Thus, by picking a generic $ P $ in $ \Geo^{(n)} , $ we can 
assume that $ \Gamma = p (P - \{v_1\}) $ consists of a chain of $ n-2 
$ great circular arcs on $ \mathcal{S} $ intersecting in four-valent 
crossings.

Suppose that $ \Gamma $ has $ c $ crossings.  Since $ \Gamma $ is 
contained in a single hemisphere of $ \mathcal{S}, $ a pair of arcs 
will intersect at most once.  Furthermore adjacent arcs cannot 
intersect, so each one of the $ n-4 $ interior arcs $ p(v_3 v_4), 
\ldots, p(v_{n-2} v_{n-1}) $ can intersect at most $ n-5 $ other arcs, 
while each of the extreme arcs $ p(v_2 v_3) $ and $ p(v_{n-1} v_n) $ 
can intersect at most $ n-4 $ other arcs.  Hence
\[
c \le \frac{1}{2}\biggl( (n-4)(n-5) + 2(n-4) \biggr) = \frac{(n-3)(n-4)}{2}.
\]

\begin{figure}[t]
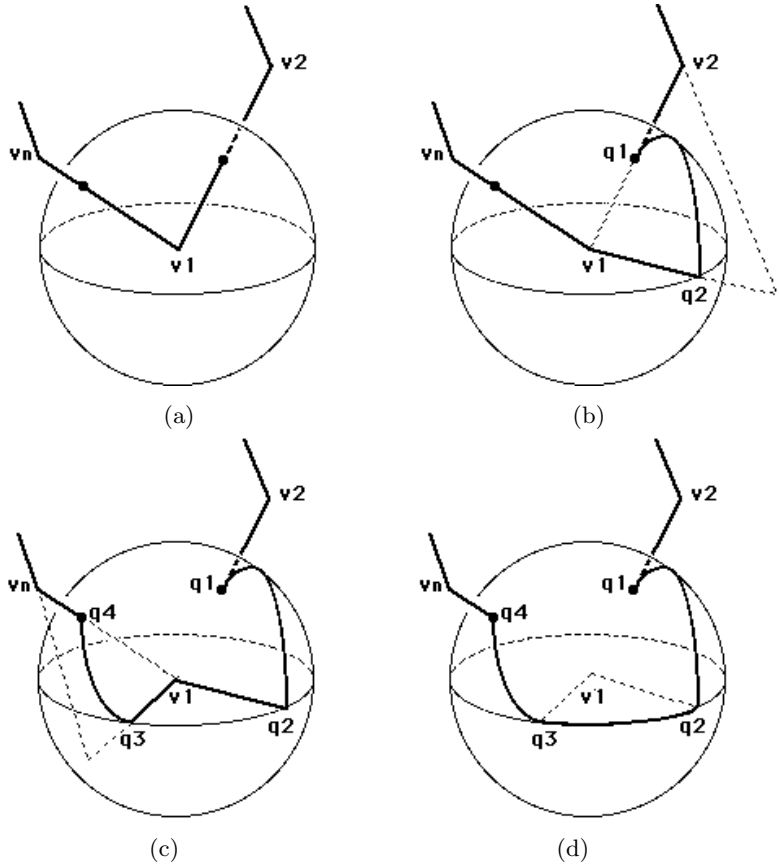

\insertfig{epsilon1.eps}{2in}
\centering{{\small (a) \hspace{48mm} (b) \hspace{4mm}}}
\insertfig{epsilon2.eps}{2in}
\centering{{\small (c) \hspace{48mm} (d) \hspace{4mm}}}
\caption{Deformation of $ P $ inside a small $ \epsilon $-ball 
about $ v_{1}.$}
\label{fig:epsilon.ball}
\end{figure}

Let $ \epsilon > 0 $ be small enough that the closed $ \epsilon $-ball 
$ \mathcal{B}_{\epsilon} $ centered at $ v_{1} $ intersects the 
polygon $ P $ in exactly two small segments of the edges $ v_{1}v_{2} 
$ and $ v_{n}v_{1} , $ as shown in Figure ~\ref{fig:epsilon.ball}(a).  
Suppose that the edge $ v_{1}v_{2} $ intersects the sphere $ \partial 
\mathcal{B}_{\epsilon} $ at the point $ q_{1} .$ Furthermore, let $ 
q_{2} $ be the point where the equator $ \partial 
\mathcal{B}_{\epsilon} \cap \mathcal{P}_{1} $ intersects the 
half-plane containing $ v_{2} $ and bounded by the line determined by 
$ v_{1} $ and $ v_{3} .$ Then we can deform the segment $ q_{1}v_{1} $ 
so that it curves along a great circle path $ \alpha_{1} $ from $ 
q_{1} $ to $ q_{2}, $ and then in a straight line path to $ v_{1} .$ 
See Figure ~\ref{fig:epsilon.ball}(b).  Note that since the arc $ 
\alpha_{1} $ lies on the same plane as $ v_{2}v_{3}, $ then $ 
p(\alpha_{1}) \cup p(v_{2}v_{3}) $ forms a single great circle 
trajectory on $ \mathcal{S} $ from $ p(v_{3}) $ to $ p(q_{2}) .$ Thus, 
after this deformation, the upper bound on the total number of 
crossings given above still holds.

Similarly, let $ q_{3} $ be the point of intersection between the 
equator $ \partial \mathcal{B}_{\epsilon} \cap \mathcal{P}_{1} $ and 
the half-plane containing $ v_{n} $ and bounded by the line determined 
by $ v_{1} $ and $ v_{n-1}, $ and let $ q_{4} $ be the point at which 
the edge $ v_{n}v_{1} $ intersects $ \partial \mathcal{B}_{\epsilon} 
.$ Then the segment $ v_{1}q_{4} $ can be deformed so that it travels 
in a straight line path from $ v_{1} $ to $ q_{3} $ and then curves 
along a great circle path $ \alpha_{3} $ from $ q_{3} $ to $ q_{4}.  $ 
See Figure ~\ref{fig:epsilon.ball}(c).  As before, the arc $ 
\alpha_{3} $ lies on the same plane as $ v_{n-1}v_{n}, $ so $ 
p(v_{n-1}v_{n}) \cup p(\alpha_{3}) $ forms a single great circle 
trajectory on $ \mathcal{S} $ from $ p(v_{n-1}) $ to $ p(q_{3}).  $ 
Therefore the upper bound on the number of crossings given above still 
holds after this deformation.

Finally, isotope $ P $ by moving $ v_1 $ into the interior of the 
triangle $ \triangle q_{2}v_{1}q_{3} $ while curving the segments $ 
q_{2}v_{1} $ and $ v_{1}q_{3} $ until they coincide with an arc along 
the equator $ \partial \mathcal{B}_{\epsilon} \cap \mathcal{P}_{1} ,$ 
as in Figure ~\ref{fig:epsilon.ball}(d).  This final transformation 
turns $ P $ into a non-polygonal embedding of the same (topological) 
knot type; this new embedding agrees with $ P $ outside of the ball $ 
\mathcal{B}_{\epsilon} $ but completely avoids its interior.  In the 
meanwhile, the image under $ p $ of this embedding is simply a 
(spherical) knot projection $ \Gamma' $ consisting of the $ n-2 $ arcs 
of $ \Gamma $ (with its ends extended by $ p(\alpha_{1}) $ and $ 
p(\alpha_{3}) $), together with an $ (n-1) $th arc $ \alpha_{2} $ 
running along the equator $ \mathcal{S} \cap \mathcal{P}_{1} $ and 
joining the endpoints $ p(q_2) $ and $ p(q_3) .$ Since $ \Gamma $ is 
contained entirely on one side of the equator, $ \alpha_{2} $ does not 
cross any other arcs.  Hence the new projection has no more crossings 
than it did before the last deformation, proving the following 
theorem.

\begin{Thm} \label{thm:new.c.bound}
Suppose that a knot $ K $ with minimal crossing number $ c(K) $ and 
minimal polygon index $ s(K) .$ Then
\begin{equation} \label{eq:crossing}
c(K) \le \frac{ \bigl(s(K) - 3 \bigr) \bigl(s(K) - 4 \bigr)}{2} . 
\end{equation}
Completing the square in (\ref{eq:crossing}) shows that 
\[ 2c \le s^2 - \, 7 s \, + \, 12 = 
\biggl(s - \frac{7}{2}\biggr)^2 - \frac{1}{4} , \]
so that 
\[ s(K) \ge \frac{7 \, + \, \sqrt{8c(K) \, + \, 1}}{2} . \]
\end{Thm}

Note that Theorem ~\ref{thm:new.c.bound} correctly predicts that the 
trefoil is the only non-trivial knot which can be realized with six 
edges.  

\begin{figure}[tp]
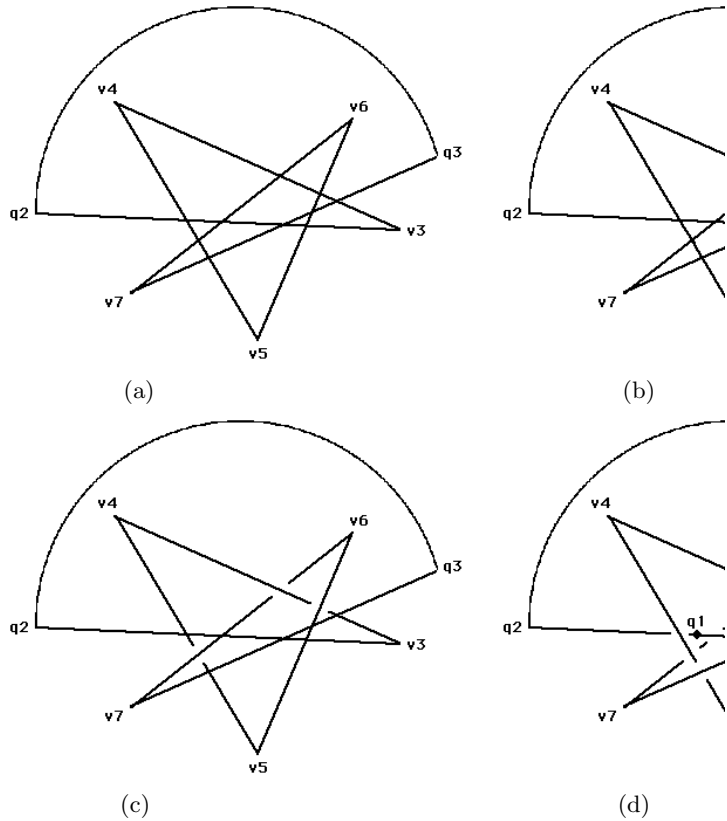

\insertfig{diagram1.eps}{1.9in}
\centering{{\small \hspace{5mm} (a) \hspace{60mm} (b) }}
\insertfig{diagram2.eps}{1.9in}
\centering{{\small \hspace{5mm} (c) \hspace{60mm} (d) }}
\caption{A knot universe and several choices in over- and 
under-crossings.}
\label{fig:diagrams}
\end{figure}

\begin{figure}[t]
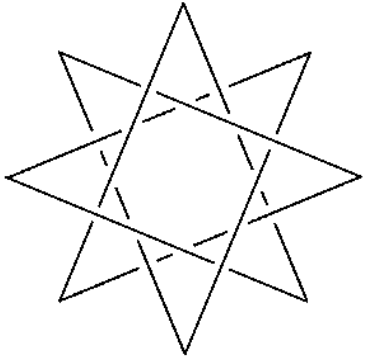

\insertfig{eight18.eps}{2in}
\caption{This octagonal embedding of the knot $ 8_{18} $ cannot be 
constructed with straight edges.}
\label{fig:eight.18}
\end{figure}

In the case of octagons, the new bound on crossing number becomes 10.  
However, in \cite{Calvo:thesis} we systematically look at the possible 
knot projections $ \Gamma' $ resulting from the deformation described 
by Figure ~\ref{fig:epsilon.ball} and thereby enumerate the 
topological knots which are appear in $ \Geo^{(8)} .$ For example, 
consider the ten-crossing knot universe shown in Figure 
~\ref{fig:diagrams}(a).  By appropriately choosing at each crossing 
which strand goes ``over'' and which one goes ``under,'' we will 
obtain a knot projection $ \Gamma' $ corresponding, as above, to some 
octagon $ P .$ As we make choices in ``over'' and ``under'' crossings 
we need to keep a few points in mind:
\begin{enumerate} \renewcommand {\labelenumi} {(\roman{enumi})} 

\item If $ v_2 v_3 $ passes under every one of its crossings, then the 
interior of the triangular disc $ \triangle v_1 v_2 v_3 $ does not 
intersect the rest of $ P .$ In this case, $ P $ can be isotoped by 
pushing $ v_2 $ in a straight line path to the midpoint of the line 
segment $ v_1 v_3 $ until $ P $ coincides with a heptagon.  A similar 
isotopy exists if $ v_{7}v_{8} $ contributes only ``under'' crossings.  
Therefore we need not consider these diagrams.

\item If the edges $ v_{5}v_{6} $ and $ v_{6}v_{7} $ both go under $ 
v_{3}v_{4}, $ as in Figure ~\ref{fig:diagrams}(b), then we can isotope 
$ P $ so that the corresponding $ \Gamma' $ has two fewer crossings.  
For instance, we can shrink the lengths of $ v_{5}v_{6} $ and $ 
v_{6}v_{7}, $ in essence performing a Reidemeister 2 move.  We can 
therefore ignore crossing choices which permit a reducing isotopy of 
this type, delaying their analysis until we examine the resulting 
reduced diagram.

\item Some choices of ``over'' and ``under'' crossings will lead to 
configurations which are impossible to create with straight edges.  
For instance, consider the three crossing choices made in Figure 
~\ref{fig:diagrams}(c).  Let $ \mathcal{P} $ be the plane containing $ 
v_{4}, v_{5}, $ and $ v_{6} .$ Note that the interior of edge $ v_6 
v_7 $ lies entirely above the plane $ \mathcal{P}, $ since it starts 
on the plane at $ v_{6} $ and then crosses over $ v_4 v_5 .$ 
Similarly, the interior of edge $ v_{3} v_{4} $ lies below the plane $ 
\mathcal{P} $ since it crosses under $ v_{5} v_{6} $ and then meets 
the plane at $ v_{4} .$ This means that $ v_{3} v_{4} $ cannot cross 
over $ v_6 v_7, $ as in Figure ~\ref{fig:diagrams}(c), unless one of 
the two edges is bent.

\item A particularly tricky example of a bad ``over'' and ``under'' 
crossing choice is shown in Figure ~\ref{fig:diagrams}(d).  This 
diagram corresponds to an octagonal realization of the knot $ 8_{18} 
,$ shown in Figure ~\ref{fig:eight.18}.  Here the problem is not as 
obvious as before.  In fact, among all of the projections 
corresponding to impossible configurations which we encounter in 
\cite{Calvo:thesis}, this is the only one which is not clearly 
impossible.  Nonetheless, through a delicate balance between 
introducing self-intersections and counting dimensions, we can show 
that there is no way to construct this configuration.  The details of 
this argument will appear in a forthcoming paper.
\end{enumerate}

After considering all possible projections $ \Gamma' $ with more than 
six crossings, we find that the only knots with polygon index 8 and 
crossing number greater than 6 are $ 8_{19} $ and $ 8_{20} .$ Since it 
is known that there are octagonal realizations of every knot $ K $ 
with crossing number $ c(K) \le 6, $ we obtain a complete list of the 
topological knots present in $ \Geo^{(8)}, $ as indicated in Theorem 
~\ref{thm:classification}(iv).  With the exception of $ 6_{3} , $ the 
square knot $ 3_{1} - 3_{1}, $ the figure-eight knot $ 4_{1}, $ and 
the unknot, every knot type in this list is chiral and therefore must 
contribute at least two path-components in $ \Geo^{(8)} .$ Therefore $ 
\Geo^{(8)} $ will contain at least twenty path-components.

\section*{Acknowledgments}

I would like to thank Ken Millett, who first led me into this 
wonderful subject, and who has always been happy to give me his 
advice, insights, and toughest questions.  I would also like to thank 
Janis Cox Millett for her hospitality this summer, as the three of us 
traveled through Paris, Athens, Delphi, Berlin, and Aix-en-Provence.

\bibliographystyle{amsplain}
\bibliography{jorge}

\end{document}

%% file: jorge.tex
\usepackage{amsfonts}

\normalsize
\theoremstyle{definition} 

\newtheorem{Thm}{Theorem}

\newcommand{\Field}[1]{\mathbb{#1}}

\newcommand{\R}{\Field{R}}
\newcommand{\Z}{\Field{Z}}

\newcommand{\norm}[1]{\lVert#1\rVert}

\input epsf.tex
  
\newcommand{\insertfig}[2]
{\begin{center} \epsfysize=#2 \epsfbox{#1} \end{center}}

%% file: hellas.bbl
\newcommand{\SortNoop}[1]{}
\providecommand{\bysame}{\leavevmode\hbox to3em{\hrulefill}\thinspace}
\begin{thebibliography}{10}

\bibitem{Adams:book}
C.~C. Adams, \emph{The knot book: An elementary introduction to the
  mathematical theory of knots}, W. H. Freeman and Co., New York, 1994.

\bibitem{Calvo:thesis}
J.~A. Calvo, \emph{Geometric knot theory: the classification of spatial
  polygons with a small number of edges}, Ph.D. thesis, University of
  California, Santa Barbara, {\SortNoop{1}}1998.

\bibitem{Calvo:hexagons}
\bysame, \emph{The embedding space of hexagonal knots}, preprint,
  {{\SortNoop{3}}1999}.

\bibitem{Calvo:montecarlo}
J.~A. Calvo and K.~C. Millett, \emph{Minimal edge piecewise linear knots},
  Ideal Knots ({A. Stasiak,} {V. Katrich,} and {L. H. Kauffman}, eds.), Series
  on Knots and Everything, vol.~19, World Scientific, Singapore,
  {\SortNoop{2}}1999, pp.~107 -- 128.

\bibitem{Adams:sticks}
{C. C. Adams,}~{B. M. Brennan,} {D. L. Greilsheimer,} and {A. K. Woo},
  \emph{Stick numbers and composition of knots and links}, Journal of Knot
  Theory and its Ramifications \textbf{6} (1997), no.~2, 149--161.

\bibitem{Furstenberg:sticks}
{E. Furstenberg,} {J. Lie,} and {J. Schneider}, \emph{Stick knots}, preprint,
  1997.

\bibitem{Kauffman:book}
L.~H. Kauffman, \emph{On knots}, Annals of Mathematics Studies, vol. 115,
  Princeton University Press, Princeton, NJ, 1987.

\bibitem{Livingston:book}
C.~Livingston, \emph{Knot theory}, Carus Mathematical Monographs, vol.~24,
  Mathematical Association of America, Washington, DC, 1993.

\bibitem{Meissen:sticks}
M.~Meissen, \emph{Edge number results for piecewise-linear knots}, Knot Theory,
  Polish Academy of Sciences, Warsaw, 1998, pp.~235--242.

\bibitem{Millett:random.knot}
K.~C. Millett, \emph{Knotting of regular polygons in 3-space}, Journal of Knot
  Theory and its Ramifications \textbf{3} (1994), no.~3, 263--278; also in
  \cite{Millett:book} pp. 31--46.

\bibitem{Millett:book}
K.~C. Millett and D.~W. Sumners (eds.), \emph{Random knotting and linking},
  Series on Knots and Everything, vol.~7, World Scientific, Singapore, 1994.

\bibitem{Negami:s.bounds}
S.~Negami, \emph{Ramsey theorems for knots, links, and spatial graphs},
  Transactions of the American Mathematical Society \textbf{324} (1991), no.~2,
  527--541.

\bibitem{Randell:conform1}
R.~Randell, \emph{A molecular conformation space}, MATH/CHEM/COMP 1987 (R.~C.
  Lacher, ed.), Studies in Physical and Theoretical Chemistry, vol.~54,
  Elsevier Science, Amsterdam, {\SortNoop{19881}}1988, pp.~125‹--140.

\bibitem{Randell:conform2}
\bysame, \emph{Conformation spaces of molecular rings}, MATH/CHEM/COMP 1987
  (R.~C. Lacher, ed.), Studies in Physical and Theoretical Chemistry, vol.~54,
  Elsevier Science, Amsterdam, {\SortNoop{19882}}1988, pp.~141‹--156.

\bibitem{Randell:sticks}
\bysame, \emph{An elementary invariant of knots}, Journal of Knot Theory and
  its Ramifications \textbf{3} (1994), no.~3, 279--286; also in
  \cite{Millett:book} pp. 47--54.

\bibitem{Rolfsen:book}
D.~Rolfsen, \emph{Knots and links}, Mathematical Lecture Series, vol.~7,
  Publish or Perish, Houston, TX, 1976.

\bibitem{Tutton:crystals}
A.~E.~H. Tutton, \emph{Crystallography and practical crystal measurement},
  vol.~1, Macmillan and Co. Ltd., London, 1922.

\bibitem{Whitney:varieties}
H.~Whitney, \emph{Elementary structure of real algebraic varieties}, Annals of
  Mathematics \textbf{66} (1967), 545‹--556.

\bibitem{Ziegler:polytopes}
G.~M. Ziegler, \emph{Lectures on polytopes}, Graduate Texts in Mathematics,
  vol. 152, Springer Verlag, New York, 1995.

\end{thebibliography}
